\documentclass[11pt, a4paper]{article}
\usepackage{amsmath,amssymb,latexsym,graphics,epsfig}
\usepackage{hyperref}
\usepackage{color}
\usepackage{amsthm}

\begin{document}
\title{A Simple and Fast Heuristic Algorithm for Edge-coloring of Graphs\footnote{Research supported by the Ministerio de Educaci\'on y
Ciencia, Spain, and the European Regional Development Fund under
project MTM2011-28800-C02-01 and by the Catalan Research Council
under project 2005SGR00256.}}

\author{M.A. Fiol and J. Vilaltella
\\ \\
{\small Universitat Polit\`ecnica de Catalunya, Dept. de Matem\`atica Aplicada IV} \\
{\small BarcelonaTech, Barcelona, Catalonia} \\
{\small (e-mails: {\tt \{fiol,joan.vilaltella\}@ma4.upc.edu})}}

\date{}

\maketitle

\begin{abstract}
\noindent A simple but empirically efficient heuristic algorithm for the edge-coloring of graphs is presented. Its basic idea is the displacement of ``conflicts" (repeated colors in the edges incident to a vertex) along paths of adjacent vertices whose incident edges are recolored by swapping alternating colors (that is, doing a Kempe interchange). The results of performance tests on random cubic and $\Delta$-regular graphs are presented, and a full implementation of the algorithm is given to facilitate its use and the reproducibility of results.
\end{abstract}

\section{Introduction}
Deciding if a graph can be $\Delta$-edge-colored (that is, whether or not it is of class one) is known to be $NP$-complete even for maximum degree $\Delta=3$ (Holyer \cite{Ho}). Therefore, assu-ming $P \ne NP$, the problem of 3-edge-coloring a cubic graph becomes intractable as the number of its vertices increases. Possible ways to  ``break" intractability are restricting the problem or admitting probabilistic solutions. Here we take the second approach, proposing a simple and fast heuristic algorithm to 3-edge-color a cubic graph $G$. Experimental data suggest that, when $G$ is 3-edge-colorable (an asymptotically almost sure property by the results of Robinson and Wormald \cite{Rob}), our algorithm efficiently gives an explicit solution and, hence, allows us to break intractability. The algorithm is easily generalized to graphs with maximum degree $\Delta$. All graphs are assumed to be connected and with no loops or multiple edges.

\section{Description of the algorithm for cubic graphs}

The heuristic algorithm is called $CVD$ (for ``conflicting vertex displacement") and starts with a random coloring of the edges of $G$ using 3 colors. A {\em conflicting vertex} is a vertex whose incident edges are not properly colored (at least two of them are equally colored). The main loop of the algorithm consists of the following steps. First, randomly choose a conflicting vertex and take one of the unproperly colored edges as the starting point of a Kempe chain. That is, a 2-edge-colored path whose two colors are the original unproper color of the starting edge and a new color that makes the vertex unconflicting or, at least, reduces the number of edges with repeated color (a suitable new color always exists, and can be randomly chosen when it is not unique). Then, swap the two colors along the chain. Such an operation does not create new conflicts but moves the conflict along the path. Stop the color swapping when, either, another conflicting vertex is reached, where the conflicts may cancel out, or the number of swaps attains a given limit (at most the number of vertices in the graph). Next, choose again a conflicting (not necessarily different) vertex and start a new Kempe chain. Repeat the procedure until, either, all conflicts are solved (so getting a 3-edge-coloring) or a fixed limit (made precise below) is reached.

The number of conflicting vertices may stabilize even if $G$ is not a snark (a class two graph, see Fiol \cite{Fiol}, Gardner \cite{Gard}, or Isaacs \cite{Is}) because not every color swap succeeds in solving one conflict. To avoid an infinite loop, we keep a counter of the number of times a new conflicting vertex is chosen and reset this counter to zero only when the number of conflicting vertices decreases. If the counter reaches a certain, constant, value $R$ or {\em repetition limit} then the main loop is restarted from a new initial random coloring of the edges. The maximum number of times the main loop can be restarted is the {\em iteration limit}, $L$. If we reach this limit, we stop and assume that the cubic graph can not be properly 3-edge-colored. Then the non-null probability of making a wrong assumption is the price for ``breaking" intractability.

\section{Complexity analysis}

Assuming the worst case, every time a conflicting vertex is chosen at most $n$ colors are swapped along a Kempe chain, and at most all $n$ vertices in the graph could be conflicting. If the number of conflicting vertices ceases to decrease, the choice of a new conflicting vertex can be repeated no more than $R$ times, the repetition limit, before the main loop restarts from a new random coloring of the edges, and the number of restarts is at most the iteration limit, $L$. The cost of randomly coloring the graph edges is $O(n)$, and at most $L$ random colorings must be done. The parameters $R$ and $L$ are fixed during the execution of the algorithm. Therefore, the expected running time is $O(n^2)+O(n)=O(n^2)$. Given the pessimistic estimation, it could be better in practice (see Section 6). A similar complexity analysis is valid in more general cases (Section 5).

\section{Algorithm modifications for $\Delta$-regular graphs and graphs with maximum degree $\Delta$}

The basic idea (conflict displacement) is the same, but this case is slightly more complex. We measure the ``conflict level'' of a vertex (how far it is from having its incident edges properly colored) by its degree ($\Delta$ if the graph is $\Delta$-regular) minus the number of different colors used in its incident edges. This gives a measure of the color repetition: assuming the colors are always in the range 0,1,2,...,$\Delta-1$, if the number of different colors used in the incident edges of a vertex is equal to the degree of that vertex, then no color is repeated and the conflict level is 0. If the same color is used in every incident edge, there are $\Delta-1$ repetitions of that color and the conflict level of the vertex is $\Delta-1$. The range of possible conflict level values is 0,1,2,...,$\Delta-1$. A data structure, a dictionary to be precise, mapping each conflict level value to the set of vertices with that conflict level, is updated as the Kempe switches change the coloring of the edges (see Section 8). This information is needed because in the new algorithm conflicting vertices are chosen according to their conflict level: vertices with the highest conflict level are chosen first. The conflict level value 0 and its corresponding set of unconflicting vertices do not need to be present in the dictionary.

Instead of the number of conflicting vertices, the total number of conflicts is used to measure the ``conflictivity'' of the graph: this is the sum of the conflict levels over all vertices in the graph, and it may be larger than the number of vertices in the graph.

To reduce the conflictivity of the initial arbitrary edge-coloring a simple greedy algorithm can be used instead of the random coloring that was sufficient in the cubic case. Besides this, and implementation details aside, the conflicting vertex choice based on the highest conflict level and the conflictivity measure based on the sum of all conflict levels are the only relevant modifications of the algorithm done to deal with the general case.

\section{Complexity analysis in the case of $\Delta$-regular graphs or graphs with maximum degree $\Delta$}

As in the previous analysis, and asuming the worst case again, every time a conflicting vertex is chosen at most $n$ swaps are done along a Kempe chain. In the present case we use the sum of all conflict levels to measure the conflictivity of the graph, which is at most $n(\Delta-1)$. If the conflictivity measure ceases to decrease, the choice of a new conflicting vertex can be repeated no more than $R$ times, the repetition limit, before the main loop restarts from a new pre-coloring of the edges, and the number of restarts is at most the iteration limit, $L$. The cost of pre-coloring all $O(\Delta n)$ graph edges using a random coloring or the greedy algorithm is $O(\Delta^2 n)$ or possibly better, depending on the efficiency of sampling from the set of $O(\Delta)$ available colors, and at most $L$ pre-colorings must be done. The parameters $R$ and $L$ are fixed during the execution of the algorithm. Therefore, the expected running time is $O((\Delta-1) n^2)+O(\Delta^2 n)=O(\Delta^2 n^2)$. As in the previous case, the estimation is pessimistic and the performance could be better in practice (see Section 7).

Roughly speaking, our approach, with its vertex-centered measure of conflict, is similar to the edge-centered approach by Lee, Wan and Guan (\cite{Lee}). These authors improved on \mbox{previous} $\Delta$-edge-coloring algorithms by giving a polynomial (albeit probabilistic) algorithm. The $CVD$ algorithm is simpler and easy to implement (see Section 8).

\section{Performance testing in the case of cubic graphs}

We tested the performance of the heuristic algorithm by measuring the time needed to 3-edge-color random cubic graphs. More specifically, for each randomly generated cubic graph, we measured the time needed to find a 3-edge-coloring, the number of iterations of the main loop, and the average time per iteration (the coloring time divided by the number of iterations). For each number of vertices, we obtained the minimum, average, and maximum values of each of these magnitudes over a set of 30 random instances, and we plotted the results. It seems that the heuristic algorithm finds a 3-edge-coloring, if any, with high probability.

\includegraphics[height=300pt,width=400pt]{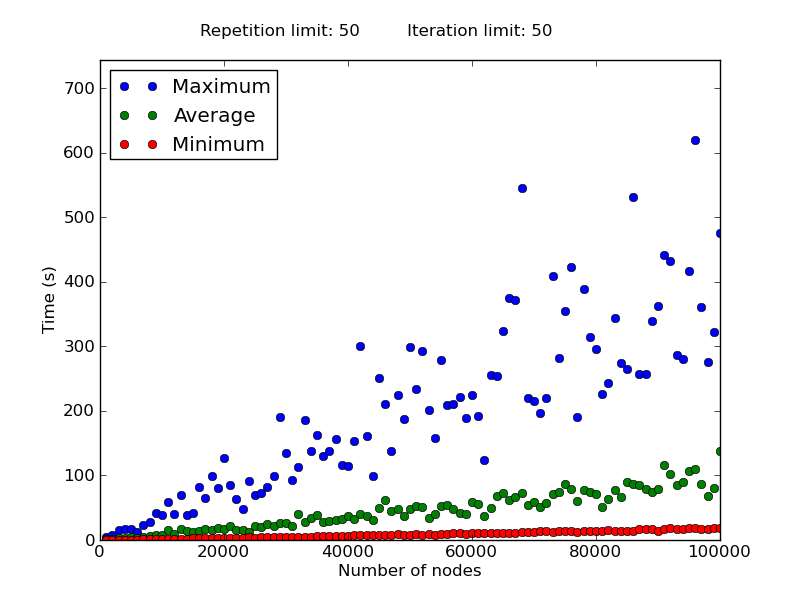}

\includegraphics[height=300pt,width=400pt]{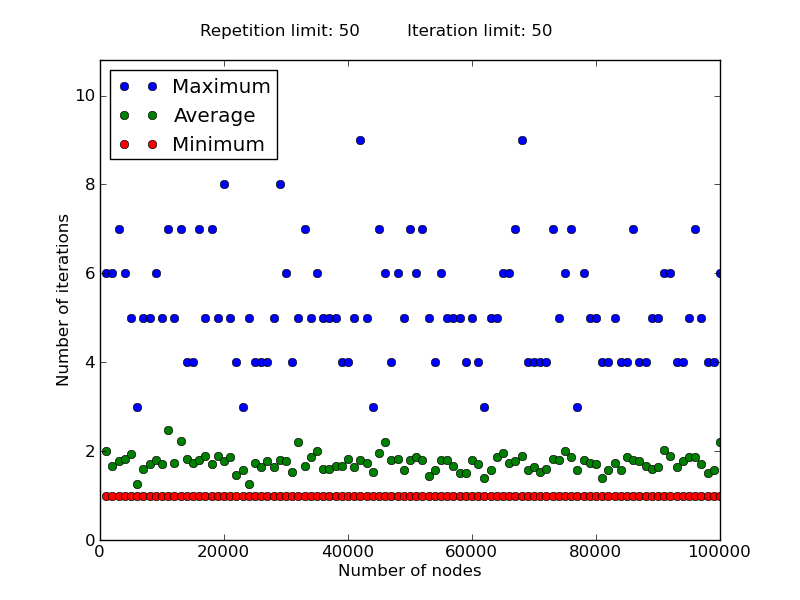}

The average running time appears to grow linearly with the number of vertices (see the first figure). The same happens with the average time per iteration (not shown), whereas the average number of iterations does not seem to depend on the number of vertices (see the second figure). A linearly growing average running time is in agreement with the $O(n^2)$ expected running time.

\section{Performance testing in the case of $\Delta$-regular graphs}

We followed the same testing procedure as in the cubic case (measuring the average running time over 30 instances), but varying the number $n$ of vertices and the value of $\Delta$. The growth of the average running time is in agreement with the expected $O(\Delta^2 n^2)$ behaviour. The following figure shows four plots, corresponding to four different values of the variable $\Delta$ (3, 7, 11 and 15).

\includegraphics[height=300pt,width=400pt]{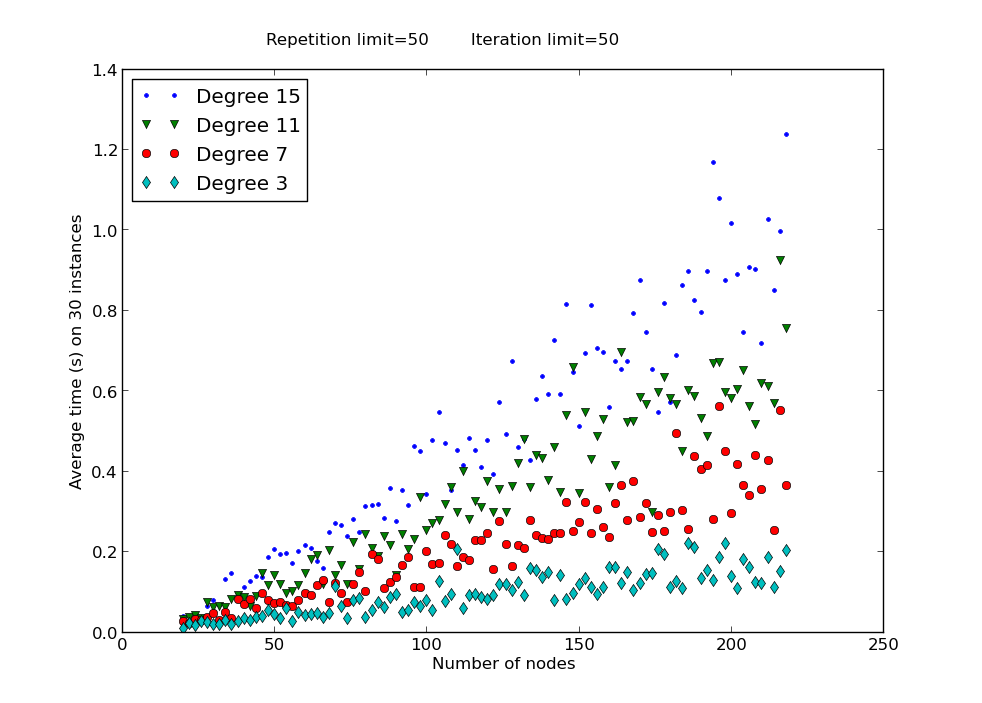}

\section{Implementation details}

We used the Python language (version 2.7.3) to implement the heuristic algorithm. This language allows us to use the NetworkX package to generate random regular graphs, and the matplotlib package to plot the results of the performance tests. 

In order to facilitate the reproducibility of our results, or the use of the CVD algorithm to anyone interested, we give a complete Python implementation of the algorithm described in Section 4, depending only on the Python Standard Library and the NetworkX package. The plot of Section 6 was created using the algorithm described in Section 3, specific for cubic graphs, with a very similar (slightly simpler) implementation.

\begin{verbatim}
### Implementation data ###
#Operating system: Linux-3.2.0-31-generic-x86_64-with-Ubuntu-12.04-precise
#CPU: Intel(R) Core(TM) i3 CPU       M 370  @ 2.40GHz. Cache size: 3072 KB
#CPU speed: 933.000 Hz

### Python code ###
from random import choice,randint
import networkx as nx

def properlyColored(G,u,D):
    return len(set(G[u].values()))==G.degree(u) and all(color in range(D) for
                                                        color in G[u].values())
def checkEdgeColoring(G,D):
    return all(properlyColored(G,u,D) for u in G.nodes())

def conflictLevel(G,u): return G.degree(u)-len(set(G[u].values()))

def createConflictDictionary(G,D):
    conflict_dictionary=dict([(i,set([])) for i in range(1,D)])
    for u in G.nodes():
        conflict_level_u=conflictLevel(G,u)
        if conflict_level_u>0: conflict_dictionary[conflict_level_u].add(u)
    return conflict_dictionary

def updateConflictDictionary(G,u,conflict_dictionary,old_conflict_level_u):
    conflict_level_u=conflictLevel(G,u)
    if old_conflict_level_u>0:
        conflict_dictionary[old_conflict_level_u].remove(u)
    if conflict_level_u>0:
        conflict_dictionary[conflict_level_u].add(u)
    return conflict_level_u-old_conflict_level_u

def maxConflictLevel(conflict_dictionary):
    return max([conflict_level for conflict_level in conflict_dictionary
                if len(conflict_dictionary[conflict_level])>0])

def totalNumberOfConflicts(conflict_dictionary):
    return sum(conflict_level*len(conflict_dictionary[conflict_level])
               for conflict_level in conflict_dictionary)

def colorEdgeAndUpdate(G,u,v,color,conflict_dictionary):
    old_conflict_level_u=conflictLevel(G,u)
    old_conflict_level_v=conflictLevel(G,v)
    G[u][v]=G[v][u]=color
    updateConflictDictionary(G,u,
                             conflict_dictionary,
                             old_conflict_level_u)
    return updateConflictDictionary(G,v,
                                    conflict_dictionary,
                                    old_conflict_level_v)


def KempeNext(G,last,node,new_color,conflict_dictionary):
    available_for_next=[w for w in G[node] if w!=last
                                           and G[node][w]==new_color]
    if available_for_next==[]: next_node=None
    else: next_node=choice(available_for_next)
    old_color=G[last][node]
    conflict_level_variation=colorEdgeAndUpdate(G,last,node,new_color,
                                                conflict_dictionary)
    return conflict_level_variation,old_color,next_node

def KempeStep(G,last,node,new_color,conflict_dictionary):
    conflict_level_variation,old_color,next_node=KempeNext(G,last,node,
                                                           new_color,
                                                           conflict_dictionary)
    if conflict_level_variation<0 or next_node==None: return node,None,None
    return node,next_node,old_color

def KempeProcess(G,last,node,new_color,conflict_dictionary):
    Kempe_chain=set([])
    while new_color!=None and last not in Kempe_chain:
        Kempe_chain.add(last)
        last,node,new_color=KempeStep(G,last,node,new_color,
                                      conflict_dictionary)

def KempeStart(G,D,node,conflict_dictionary):
    colors=set(range(D))
    next_node=None
    for adjacent in G[node]:
        edge_color=G[node][adjacent]
        if edge_color in colors: colors.remove(edge_color)
        else: next_node=adjacent
    if next_node!=None:
        KempeProcess(G,node,next_node,choice(list(colors)),conflict_dictionary)

def preColoring(G,D): #Pre-coloring with a greedy algorithm
    for e in G.edges(): G[e[0]][e[1]]=G[e[1]][e[0]]=None
    for e in G.edges():
        available_colors=set(range(D))
        available_colors-=set(G[e[0]].values())
        available_colors-=set(G[e[1]].values())
        if available_colors==set(): G[e[0]][e[1]]=G[e[1]][e[0]]=randint(0,D-1)
        else: G[e[0]][e[1]]=G[e[1]][e[0]]=choice(list(available_colors))


#def preColoring(G,D): #Random pre-coloring
#   for e in G.edges(): G[e[0]][e[1]]=G[e[1]][e[0]]=randint(0,D-1)

def heuristic(G,D,repetition_limit):
        repetitionCounter=0
        conflict_dictionary=createConflictDictionary(G,D)
        previous=current=totalNumberOfConflicts(conflict_dictionary)
        while previous>0:
            highest_conflict_level=maxConflictLevel(conflict_dictionary)
            node=choice(list(conflict_dictionary[highest_conflict_level]))
            KempeStart(G,D,node,conflict_dictionary)
            current=totalNumberOfConflicts(conflict_dictionary)
            if current==0: return True
            if current>=previous:
                repetitionCounter+=1
                if repetitionCounter>repetition_limit: return False
            else: repetitionCounter=0
            previous=min(previous,current)
        return True

def applyHeuristic(G,D,repetition_limit,iteration_limit):
    preColoring(G,D)
    number_of_iterations=1
    while not heuristic(G,D,repetition_limit):
        if number_of_iterations>iteration_limit: break
        preColoring(G,D)
        number_of_iterations+=1
    print "Number of iterations:",number_of_iterations
    print "Edge-coloring successful:",checkEdgeColoring(G,D)

### Example ###
repetition_limit=iteration_limit=50
NUMBER_OF_VERTICES=50000
DEGREE=3
G=nx.random_regular_graph(DEGREE,NUMBER_OF_VERTICES)
applyHeuristic(G,DEGREE,repetition_limit,iteration_limit)

\end{verbatim}

\section{Testing a conjecture on odd graphs}

Odd graphs, $O_k$ with $k \geq 2$ an integer, are defined in the following way: the vertices correspond to the $(k-1)$-subsets of a $(2k-1)$-set, and two vertices are adjacent if their corresponding subsets are disjoint. It is conjectured (Biggs \cite{Biggs}, Fiorini and Wilson \cite{Fior}) that odd graphs are of class one except for $k=3$ (the Petersen graph) and $k$ a power of two (when the graph has an odd number of vertices, implying that it is trivially of class two). In their book, Fiorini and Wilson comment results for $k \leq 8$. We used the $CVD$ algorithm to test that the conjecture is true for $k \leq 11$.

\section{Conclusions and future directions}

We have seen that random cubic and $\Delta$-regular graphs can be edge-colored with empirical efficiency by a relatively simple heuristic algorithm (with simple data structures).

Many combinatorial problems, for instance logical circuit problems (Fiol \cite{Fiol}), can be \mbox{reduced} to 3-edge-coloring of cubic graphs or edge-coloring of regular graphs. Therefore, any efficient edge-coloring heuristic might be useful to solve them. Moreover, efficient edge-coloring is of practical interest in many applications.

We are currently working in the generalization of the heuristic algorithm to hypergraphs and other combinatorial structures, and testing its performance in combinatorial problems.


\end{document}